\documentclass{amsart}

\usepackage{amsmath,amssymb,graphicx}

\newtheorem{theorem}{Theorem}
\newcommand{\bt}{\begin{theorem}}
\newcommand{\et}{\end{theorem}}

\newtheorem*{theoremNN}{Theorem}
\newcommand{\btNN}{\begin{theoremNN}}
\newcommand{\etNN}{\end{theoremNN}}

\newtheorem{lemma}{Lemma}
\newcommand{\bl}{\begin{lemma}}
\newcommand{\el}{\end{lemma}}

\newtheorem{corollary}{Corollary}
\newcommand{\bc}{\begin{corollary}}
\newcommand{\ec}{\end{corollary}}

\newtheorem{definition}{Definition}
\newcommand{\bdf}{\begin{definition}}
\newcommand{\edf}{\end{definition}}

\newtheorem{conjecture}{Conjecture}
\newcommand{\bconj}{\begin{conjecture}}
\newcommand{\econj}{\end{conjecture}}

\newtheorem*{conjectureNN}{Conjecture}
\newcommand{\bconjNN}{\begin{conjectureNN}}
\newcommand{\econjNN}{\end{conjectureNN}}

\newtheorem{example}{Example}
\newcommand{\bex}{\begin{example}}
\newcommand{\eex}{\end{example}}

\newtheorem{problem}{Problem}
\newcommand{\bprob}{\begin{problem}}
\newcommand{\eprob}{\end{problem}}

\newtheorem{oproblem}{Open Problem}
\newcommand{\boprob}{\begin{oproblem}}
\newcommand{\eoprob}{\end{oproblem}}

\newcommand{\beq}{\begin{equation}}
\newcommand{\eeq}{\end{equation}}
\newcommand{\benum}{\begin{enumerate}}
\newcommand{\eenum}{\end{enumerate}}

\newcommand{\N}{\ensuremath{ \mathbf N }}

\newcommand{\R}{\ensuremath{\mathbf R}}

\newcommand{\Rn}{\ensuremath{ \mathbf{R}^n }}

\newcommand{\Z}{\ensuremath{\mathbf Z}}

\newcommand{\Zn}{\ensuremath{ \mathbf{Z}^n }}

\newcommand{\mcg}{\ensuremath{ \mathcal G}}

\newcommand{\mcs}{\ensuremath{ \mathcal S}}

\newcommand{\bq}{\begin{eqnarray*}}
\newcommand{\eq}{\end{eqnarray*}}
\newcommand{\be}{\begin{eqnarray}}
\newcommand{\ee}{\end{eqnarray}}
\newcommand{\ba}{\begin{array}}
\newcommand{\ea}{\end{array}}
\newcommand{\bfr}{\begin{flushright}}
\newcommand{\efr}{\end{flushright}}

\newcommand{\bmat}{\left(\begin{matrix}}
\newcommand{\emat}{\end{matrix}\right)}

\DeclareMathOperator{\conv}{\text{conv}}

\DeclareMathOperator{\qqand}{\qquad\text{and}\qquad}

\title[Asymptotic approximate groups]
{Every finite subset of an abelian group is an asymptotic approximate group}
\author{Melvyn B. Nathanson}
\address{Department of Mathematics\\Lehman College (CUNY)\\Bronx, NY 10468} \email{melvyn.nathanson@lehman.cuny.edu}

\thanks{Supported in part by a grant from the PSC-CUNY Research Award Program.}

\date{\today}

\begin{document}

\maketitle

\begin{abstract}
If $A$ is a nonempty subset of an additive abelian group $G$, 
then the \emph{$h$-fold sumset} is 
\[
hA = \{x_1  + \cdots + x_h : x_i \in A_i \text{ for } i=1,2,\ldots, h\}.
\]  
We do not assume that $A$ contains the identity, nor that $A$ is symmetric, 
nor that $A$ is finite.
The set $A$  is an \emph{$(r,\ell)$-approximate group in $G$} 
if there exists a subset $X$ 
of $G$ such that $|X| \leq \ell$ and $rA \subseteq XA$.
The set $A$ is an \emph{asymptotic $(r,\ell)$-approximate group} 
if the sumset $hA$ is an $(r,\ell)$-approximate group for all sufficiently large $h$.  
It is proved that every polytope in a real vector space is an asymptotic 
$(r,\ell)$-approximate group, that every finite set of lattice points 
is an asymptotic $(r,\ell)$-approximate group, 
and that every finite subset of every abelian group is 
an asymptotic $(r,\ell)$-approximate group.
\end{abstract}

\section{Approximate groups}
Let $G$ be a group, not necessarily abelian, and written multiplicatively.   
Let $\N_0$ and \N\ denote the sets of nonnegative integers 
and positive integers, respectively.  
For  subsets  $A_1,\ldots, A_h$, and $A$ of $G$, we define the \emph{product sets}
\[
A_1\cdots A_h = \{a_1\cdots a_h: a_i \in A_i \text{ for } i=1,\ldots, h\} 
\]
and
\[
A^h = \underbrace{A\cdots A}_{\text{$h$ factors}} = \{a_1\cdots a_h: a_i \in A \text{ for } i=1,\ldots, h\}.
\]
We have 
\[
\left(A^h\right)^r = A^{hr} \qqand A^{h+r} = A^h A^r 
\]
for all positive integers $h$ and $r$.

For $c \in G$, we define the  \emph{left} and  \emph{right translates} 
\[
c A = \{ca:a \in A \}
\qqand
A  c = \{ac:a \in A \}.  
\]
The \emph{normalizer} of $A$ is the subgroup 
\[
N_G(A) = \{c\in G: c  A = A  c\}.
\]  
If $c\in N_G(A)$ and $r \in \N$, then 
\[
(c  A)^r = c^r   A^r.
\]
If $G$ is abelian, then for all $c \in G$ and $A \subseteq G$ 
we have $c A = A  c$, and so $N_G(A) = G$.

Let $A$ be a nonempty subset of $G$, and let $r$ and $\ell$ be positive integers.  
The set $A$ is an \emph{$(r,\ell)$-approximate group} 
if there exists a set $X \subseteq G$ such that 
\beq     \label{AAG-def1}
|X| \leq \ell
\eeq
and
\beq     \label{AAG-def2}
A^r \subseteq XA.
\eeq
The idea of an approximate group evolved from Freiman's inverse theorem 
in additive number theory (Freiman~\cite{frei64,frei64t,frei66,frei73}, 
Nathanson~\cite{nath96bb}).
The definition of approximate group  in this paper is less restrictive  
than the original definition, which appears in Tao~\cite{tao07} 
and which has been extensively investigated (for example, by 
Breuillard, Green, and Tao~\cite{breu13,breu15,breu-gree-tao12,gree12,gree14} 
and Pyber and Szabo~\cite{pybe-szab14}).  
Helfgott~\cite{helf15} is a recent survey of related problems on growth in groups.

The set $A$ is an \emph{asymptotic $(r,\ell)$-approximate group} 
if every sufficiently high power of $A$ is an $(r,\ell)$-approximate group.
This means that there exists an integer $h_0(A)$ such that, for each integer
 $h \geq h_0(A)$, 
there is a set $X_h \subseteq G$ such that 
\beq     \label{AAG-def1-asymp}
|X_h| \leq \ell
\eeq
and
\beq     \label{AAG-def2-asymp}
A^{hr} \subseteq X_h A^h.
\eeq

We begin with some simple facts about approximate groups.

\bl              \label{AAG:lemma:translate}
Let $A$ be a nonempty subset of a group $G$, and let $c \in N_G(A)$.
Let $r$ and $\ell$ be positive integers.  
If $A$ is an $(r,\ell)$-approximate group, 
then $cA$ is an $(r,\ell)$-approximate group.  
If $A$ is an asymptotic $(r,\ell)$-approximate group, 
then $cA$ is an asymptotic $(r,\ell)$-approximate group.  
\el

\begin{proof}
Suppose that $h \geq 1$ and that $X_h$ is a subset of $G$ 
that satisfies~\eqref{AAG-def1-asymp} 
and~\eqref{AAG-def2-asymp}.  Let 
\[
X'_h = c^{rh} X_h c^{-h}.
\]
We have  
\[
|X'_h| = |X_h| \leq \ell
\]
and
\[
(cA)^{rh} = c^{rh}A^{rh} \subseteq c^{rh}X_h A^h 
= \left( c^{rh} X_h c^{-h} \right) c^hA^h = X'_h (cA)^h.
\]
This completes the proof.  
\end{proof}

\bl        \label{AAG:lemma:DirectProduct}
If $A_0$ is a $(r,\ell_0)$-approximate subgroup of the group $G_0$ 
and $A_1$ is a $(r,\ell_1)$-approximate subgroup of $G_1$, 
then $A_0 \times A_1$ is a $(r,\ell_0\ell_1)$-approximate subgroup 
of $G_0 \times G_1$.
\el

\begin{proof}
There exist subsets $X_0$ of $G_0$ and $X_1$ of $G_1$ 
such that $|X_0| \leq \ell_0$, $|X_1| \leq \ell_1$, and 
\[
A_0^r \subseteq X_0A_0 \qqand A_1^r \subseteq X_1A_1.
\]
It follows that 
\begin{align*}
(A_0 \times A_1)^r 
& = A_0^r \times A_1^r \\
&\subseteq X_0A_0 \times X_1A_1 \\
& = (X_0\times X_1)(A_0\times A_1)  
\end{align*}
and $| X_0\times X_1 | = | X_0| |X_1 | \leq \ell_1\ell_1$.  
This completes the proof.  
\end{proof}

\bl         \label{AAG:lemma:Main-projection}
Let $G_0$ and $G_1$ be groups such that $G_0$ is finite  with $|G_0| = n_0$.
Let  $\pi_1: G_0 \times G_1 \rightarrow G_1$ 
be the projection homomorphism 
defined by $\pi_1(u_0,u_1) = u_1$ for all $(u_0,u_1) \in G_0 \times G_1$.
Let $A$ be a  subset of $G_0 \times G_1$.  
If $\pi_1(A)$ is an $(r,\ell)$-approximate group in $G_1$, 
then $A$ is an $(r, n_0\ell) $-approximate group in $G_0 \times G_1$.
If $\pi_1(A)$ is an asymptotic $(r,\ell)$-approximate group in $G_1$, 
then $A$ is an asymptotic $(r,\ell)$-approximate group in $G_0 \times G_1$.
\el

\begin{proof}
Let $A$ be a nonempty subset of $G$, and let 
$A_1 = \pi_1(A) \subseteq G_1$.  
If $A_1$ is an $(r,\ell)$-approximate group in $G_1$, 
then there exists $X_1 \subseteq G_1$ such that 
\[
|X_1| \leq \ell
\]
and 
\beq                                                 \label{AAG:Main-projection}
A_1^r \subseteq X_1  A_1.
\eeq
Let 
\[
X = G_0\times X_1.
\]
We have
\[
|X| = |G_0| |X_1| \leq n_0\ell.
\]
If $(b_0,b_1) \in A^r$, then $b_0 \in  G_0$ and $b_1 \in A_1^r$.  
It follows from~\eqref{AAG:Main-projection} 
that there exist $x_1 \in X_1$ and $a_1 \in A_1$ such that 
\[
b_1 = x_1a_1.
\]
Because $A_1 = \pi_1(A)$, there exists $a_0 \in G_0$ such that $(a_0,a_1) \in A$.  
Then 
\[
(b_0,b_1)  = (b_0, x_1a_1) = (b_0 a_0^{-1},x_1)( a_0, a_1) \in (G_0\times X_1) A = XA
\]
and so 
\[
A^r \subseteq XA.
\]
This proves that $A$ is an $(r, n_0\ell) $-approximate group in $G$.

If $A_1$ is an asymptotic $(r,\ell)$-approximate group in $G_1$,  
then there is an integer $h_0(A_1)$ such that, for all $h \geq h_0(A_1)$, the set 
\[
hA_1 = h\pi_1(A) = \pi_1(hA)
\]
is an $(r,\ell)$-approximate group in $G_1$, 
and so $hA$ is an $(r, n_0 \ell)$-approximate group in $G_0 \times G_1$.  
Thus, $A$ is an asymptotic $(r, n_0\ell) $-approximate group in $G$.  
This completes the proof.  
\end{proof}

\bl         \label{AAG:lemma:FiniteGroup}
Let $G_0$ be a finite group with $|G_0| = n_0$, and let $r \in \N$.  
Every nonempty subset of $G_0$ is an $(r,n_0)$-approximate group.
\el

\begin{proof}
Let $X = G_0$.  If $A$ is a nonempty subset of $G_0$, then 
\[
A^r \subseteq G_0 = XA 
\]
and $|X| = |G_0| = n_0$.  This completes the proof.  
\end{proof}

\bl             \label{AAG:lemma:isomorphism}
Let $G$ and $G'$ be groups, and let $f:G\rightarrow G'$ be a homomorphism.  
Let $A$ be a nonempty subset of $G$.   
If $A$ is an $(r,\ell)$-approximate group in $G$,   
then $f(A)$ is an $(r,\ell)$-approximate group in $G'$.  
If $A$ is an asymptotic $(r,\ell)$-approximate group in $G$,   
then $f(A)$ is an asymptotic $(r,\ell)$-approximate group in $G'$.  
\el

\begin{proof}
It suffices to observe that if $f:G\rightarrow G'$ is a homomorphism, 
and if $A$ and $X$ are subsets of $G$ such that $|X| \leq \ell$ 
and $A^r \subseteq XA$ for some $r \geq 2$,  
then $|f(X)| \leq \ell$ and 
\[
\left( f(A) \right)^r = f\left( A^r \right) \subseteq f(XA) = f(X) f(A).
\]
This completes the proof. 
\end{proof}

The following example  shows that there are finite subsets of groups 
that are not  asymptotic $(r,\ell)$-approximate groups for any 
integers $r \geq 2$ and $\ell \geq 1$.    

\bt
There is a nonabelian group $G$ and a nonempty finite subset 
$A$ of $G$ such that, for all integers $r \geq 2$ and $\ell \geq 1$, 
the set $A$ is not an  asymptotic $(r,\ell)$-approximate group.
\et

\begin{proof}
Let $G$ be the free group of rank 2 generated by the set $A = \{a_0, a_1\}$.
For all $n \in \N$ we have $|A^n| = 2^n$.  
Let $h, \ell, r \in \N$ with $r \geq 2$ and 
\[
h > \frac{\log_2 \ell}{r-1}.
\]
If $X_h$ is a finite subset of $G$ such that 
\[
A^{hr} \subseteq X_h A^h
\]  
then
\[
2^{hr} = \left| A^{hr} \right|  \leq   \left|  X_h A^{h} \right|  
\leq  |X_h|  \left| A^{h} \right| = |X_h| 2^h
\]
and so 
\[
|X_h| \geq 2^{h(r-1)} > \ell.
\]
Thus, $A$ is not an asymptotic $(r,\ell)$-approximate group.
\end{proof}

In this paper we study abelian groups, written additively. 
For subsets $A_1,\ldots, A_h$, and $A$ of an additive abelian group $G$, 
we define the \emph{sumsets}  
\[
A_1+ \cdots + A_h = \{a_1 + \cdots + a_h: a_i \in A_i \text{ for } i=1,\ldots, h\}
\]
and
\[
hA = \underbrace{A + \cdots + A}_{\text{$h$ summands}} 
= \{a_1 + \cdots + a_h: a_i \in A \text{ for } i=1,\ldots, h\}.
\]
Let  $\mcs(A) = \bigcup_{h=1}^{\infty} hA$
be the subsemigroup of $G$ generated by $A$,  
and let $\mcg(A)$ be  the subgroup of $G$ generated by $A$. 
The set $A$ is an \emph{$(r,\ell)$-approximate group} 
if there exists a set $X \subseteq G$ such that 
\beq     \label{AAG-def1abel}
|X| \leq \ell
\eeq
and
\beq     \label{AAG-def2abel}
rA \subseteq X + A.
\eeq
The set $A$ is an asymptotic \emph{$(r,\ell)$-approximate group} 
if there exists an integer $h_0(A)$ such that $hA$ is an 
$(r,\ell)$-approximate group for every $h \geq h_0(A)$.  
This means that,  for every $h \geq h_0(A)$, 
there exists a set $X_h \subseteq G$ such that 
\[
|X_h| \leq \ell
\]
and
\[
rhA \subseteq X_h + hA.
\]
Nathanson~\cite{nath16c} proved that every finite set of integers 
is an asymptotic approximate group.  
The goal of this paper is to prove that every  finite subset 
of every abelian group is an asymptotic approximate group.

\section{Polytopes and approximate groups}

Let $V$ be a real vector space.  The \emph{dilation} of a subset $X$ of $V$
by the real number $\lambda$ is the set
\[
\lambda \ast X = \{\lambda x:x \in X\}.
\]

\bl                       \label{approx:lemma:ConvexSumset}
Let $K$ be a convex subset of a real vector space $V$.  
For every positive integer $h$, 
\[
h\ast K = hK.
\]
\el

Thus, the $h$-fold sumset of a convex set is the dilation of the set by $h$.

\begin{proof}
If $z \in h\ast K$, then there exists $x \in K$ such that $z = hx$.
Because
\[
hx = \underbrace{x + x + \cdots + x}_{\text{$h$ summands}}
\]
it follows that $z \in hK$, and so $h \ast K \subseteq hK$.

Conversely, if $z \in hK$, then there exist $x_1,\ldots, x_h \in K$ such that 
\[
z = x_1 + \cdots + x_h.
\] 
Because the set $K$ is convex, it contains the convex combination 
\[
y = \frac{1}{h}x_1 + \cdots + \frac{1}{h}x_h
\]
and so $z = hy \in h\ast K$.  Therefore, $hK \subseteq h\ast K$.
This completes the proof.  
\end{proof}

In a real vector space $V$, a \emph{polytope} is the convex hull 
of a nonempty  finite set of points.  
We denote by $\Z^n$ the group of lattice points in the vector space \Rn.
A \emph{lattice polytope} in \Rn\ is the convex hull 
of a nonempty  finite set of  lattice points.  

An additive subgroup $\mcg$ of a real normed vector space $V$ is \emph{discrete}
if every bounded subset of $V$ contains only finitely many elements of $\mcg$.
For example, \Zn\ is a discrete subgroup of \Rn.
Every subgroup of a discrete group is discrete.  
The subgroup of \R\ generated by the set $\{1, \sqrt{2} \}$ is not discrete.

For integers $r \geq 2$ and $k \geq 2$, define the binomial coefficient
\[ 
b(r,k) = \binom{(r+1)(k-1)-1}{k-1}.
\]
Note that $b(r,2) = r$ for all $r \geq 2$.

\bt                   \label{AAG:theorem:PolytopeApproxGroup}
Let  $A$ be a finite subset of a  real vector space $V$, 
with $|A| = k \geq 2$.  
Let $\mcs(A)$ be the subsemigroup of $V$ generated by $A$.
For every integer $r \geq 2$, the polytope 
\[
K = \conv(A)
\]
is an $(r,b(r,k))$-approximate group in $V$.  
Moreover, there exists a set $X_r \subseteq \mcs(A)$ such that 
\beq     \label{AAG:TightBound}
|X_r| \leq b(r,k) 
\eeq
and
\beq                \label{AAG:PolytopeApproxGroup}
rK \subseteq \frac{1}{k-1} \ast X_r+K.
\eeq
\et

\begin{proof}
Let $A = \{a_0, a_1,\ldots, a_{k-1} \}$.  We have  
\begin{align*}
K & = \conv(A) 
 =  \left\{  \sum_{i=0}^{k-1} \lambda_i a_i: \lambda_i \geq 0 
\text{ and } \sum_{i=0}^{k-1} \lambda_i =1  \right\} \\
& =  a_0 +  \left\{  \sum_{i=1}^{k-1} \lambda_i (a_i-a_0): 
\lambda_i \geq 0 
\text{ and } \sum_{i=1}^{k-1} \lambda_i \leq 1  \right\} \\
& = a_0 + K_0 
\end{align*}
where 
\[
A_0 = A- a_0 = \{0\} \cup \{a_i - a_0:i=1,\ldots, k-1\} 
\]
and 
\begin{align*}
K_0 & = \conv(A_0) = \left\{  \sum_{i=1}^{k-1} \lambda_i (a_i-a_0): 
\lambda_i \geq 0 \text{ and } \sum_{i=1}^{k-1} \lambda_i \leq 1  \right\} . 
\end{align*}
By Lemma~\ref{AAG:lemma:translate}, the polytope $K_0$ is an 
$(r, \ell )$-approximate group if and only if $K$ is an 
$(r, \ell )$-approximate group.  
Moreover, if $X'_r \subseteq \mcs(A)$ and 
\[
rK_0 \subseteq \frac{1}{k-1} \ast X'_r+K_0
\]
then
\begin{align*}
rK & = r(a_0+K_0) = ra_0 + rK_0  \\
&  \subseteq ra_0 +  \frac{1}{k-1} \ast X'_r+K_0 \\
&  = (r-1)a_0 +  \frac{1}{k-1} \ast X'_r+ a_0 + K_0 \\
& =  \frac{1}{k-1} \ast X_r + K
\end{align*}
where 
\[
X_r =  (k-1)(r-1)a_0 + X'_r 
\]
 is a subset of $\mcs(A)$ and 
$|X_r | = |X'_r|$.
Thus, we can assume that $A = \{0, a_1,\ldots, a_{k-1} \}$, and so 
\[
K = \left\{ \sum_{i=1}^{k-1} \lambda_i a_i: \lambda_i \geq 0  \text{ and }
\sum_{i=1}^{k-1} \lambda_i \leq 1 \right\} 
\]
and 
\[
rK = \left\{ \sum_{i=1}^{k-1} \mu_i a_i: \mu_i \geq 0\text{ and }
\sum_{i=1}^{k-1} \mu_i \leq r \right\}.
\]

Consider the finite set 
\begin{align*}
X_r & = \left\{ \sum_{i=1}^{k-1} m_i  a_i :(m_1,\ldots, m_{k-1}) \in \N_0^{k-1} 
\text{ and } \sum_{i=1}^{k-1} m_i \leq r(k-1) - 1\right\}.
\end{align*}
Note that  $X_r \subseteq \mcs(A)$.  
The number of  $(k-1)$-tuples $(m_1,\ldots, m_{k-1}) \in \N_0^{k-1}$ 
such that $\sum_{i=1}^{k-1} m_i \leq r(k-1) - 1$ is 
the binomial coefficient $ b(r,k) = \binom{(r+1)(k-1)-1}{k-1}$, and so 
\[
|X_r| \leq b(r,k).  
\]

Let 
\[
w =  \sum_{i=1}^{k-1} \mu_i a_i \in rK
\]
with $\mu_i \geq 0$ for $i=1,\ldots, k-1$ and $\sum_{i=1}^{k-1} \mu_i \leq r$.   
For each $i \in \{1,2,\ldots, k-1\}$ 
there is a nonnegative integer $m_i$ 
and a nonnegative real number $\lambda_i$ such that 
\[
\frac{m_i}{k-1} \leq \mu_i < \frac{m_i+1}{k-1}
\]
and
\[
0 \leq \lambda_i = \mu_i -  \frac{m_i}{k-1} < \frac{1}{k-1}. 
\]
We have 
\[
\sum_{i=1}^{k-1} \lambda_i < 1
\]
and so  
\[
\sum_{i=1}^{k-1} \lambda_i a_i \in K.
\]
Similarly,  
\[
\frac{1}{k-1}  \sum_{i=1}^{k-1} m_i  
\leq  \frac{1}{k-1}  \sum_{i=1}^{k-1} m_i + \sum_{i=1}^k \lambda_i 
= \sum_{i=1}^{k-1} \mu_i \leq r
\]
and so
\beq   \label{AAG:mrkIneq}
\sum_{i=1}^{k-1} m_i \leq r(k-1)
\eeq

There are two cases.  If
\[
\sum_{i=1}^{k-1} \lambda_i > 0
\]
then 
\[
 \frac{1}{k-1}  \sum_{i=1}^{k-1} m_i <  r
\]
and so 
\[
\sum_{i=1}^{k-1} m_i <  r(k-1).
\]
Because  $m_1,\ldots, m_{k-1}$ are nonnegative integers, we have 
\[
\sum_{i=1}^{k-1} m_i \leq r(k-1) -1
\]
and 
\begin{align*}
w  &  =  \sum_{i=1}^{k-1} \mu_i a_i  
 = \sum_{i=1}^{k-1}  \left(  \frac{m_i}{k-1} + \lambda_i \right)  a_i \\ 
& =  \frac{1}{k-1}  \sum_{i=1}^{k-1} m_i a_i + \sum_{i=1}^k \lambda_i a_i 
\in  \frac{1}{k-1} \ast X_r+K.
\end{align*}

In the second case, 
\[
\sum_{i=1}^{k-1} \lambda_i = 0  
\]
and so 
\[
w =   \frac{1}{k-1} \sum_{i=1}^{k-1} m_i a_i.   
\]
Note that $0 \in  \left(1/(k-1)\right) \ast X_r+K$.
If $w \neq 0$, then $m_j \geq 1$ for some $j \in \{1,\ldots, k-1\}$.
It follows from inequality~\eqref{AAG:mrkIneq} that 
\[
\sum_{\substack{i=1\\i \neq j}}^{k-1} m_i   + (m_j-1)
\leq r(k-1) -1
\]
and so 
\[
\sum_{\substack{i=1\\i \neq j}}^{k-1} m_i a_i  + (m_j-1)a_j \in X_r.  
\]
Moreover, 
\[
\frac{1}{k-1} a_j \in K.
\]
Therefore, 
\[
w =   \frac{1}{k-1} \left( \sum_{\substack{i=1\\i \neq j}}^{k-1} m_i a_i  
+ (m_j - 1)a_j \right) +  \frac{1}{k-1} a_j  
\in  \frac{1}{k-1} \ast X_r + K.
\]
This completes the proof.  
\end{proof}

The following example shows that the upper bound~\eqref{AAG:TightBound} is tight.  
Let $k=2$ and let $A = \{ a_0, a_1\} \in \R$ with $a_0 < a_1$.  Then  
\[
K = \conv(A) = \{(1-t)a_0+ta_1: 0 \leq t \leq 1\} = [ a_0, a_1 ].  
\]
For $r \geq 2$, we partition the sumset $rK$ as follows:
\begin{align*}
rK & = [r a_0, r a_1]   \\
& =   \bigcup_{i=1}^r [  (r - i + 1) a_0 + (i-1) a_1,  (r - i) a_0 + ia_1 ]   \\
& =   \bigcup_{i=1}^r  \left( (r - i ) a_0 + (i-1) a_1 +[a_0,a_1]  \right) \\
& = X_r + K
\end{align*}
where 
\[
X_r=  \{  (r-i)a_0 + (i-1)a_1 : i=1,\ldots, r \} 
\]
and
\[
  |X_r| = r = b(r,k).  
 \]
 Conversely, let $X$ be any finite subset of $\R$  such that $rK \subseteq X+K$.  
 Denoting the Lebesgue measure  of a set $S$ by $\mu(S)$, we obtain 
 \[
 \mu(K) = a_1-a_0  > 0
 \]
 and
 \[
 r\mu(K)= \mu(rK) \leq \mu(X + K) \leq \sum_{x \in X} \mu(x + K) = |X| \mu(K)
 \]
 and so $|X| \geq r$.  Thus,~\eqref{AAG:TightBound} is tight.

\bt                   \label{AAG:theorem:PolytopeAsymptoticApproxGroup}
Let  $A$ be a finite subset of a  real vector space $V$, 
with $|A| = k \geq 2$.  
Let $\mcs(A)$ be the subsemigroup of $V$ generated by $A$. 
Consider the polytopes $K = \conv(A)$ and $K' = \conv( (k-1) \ast A) = (k-1)K$.
For every pair $(r,h)$ of integers with $r \geq 2$ and $h \geq 1$,
there exists a set 
\[
X_{r,h} \subseteq \mcs(h\ast A) \subseteq \mcs(A)
\]
such that 
\beq           \label{AAG:PolytopeAsymptoticApproxGroup-X}
|X_{r,h} | \leq b(r,k) 
\eeq
\beq           \label{AAG:PolytopeAsymptoticApproxGroup}
rhK \subseteq \frac{1}{k-1} \ast X_{r,h}+ hK
\eeq
and
\beq           \label{AAG:PolytopeAsymptoticApproxGroup-K'}
r h K' \subseteq X_{r,h}  + h K'.
\eeq
In particular, $K$ and $K'$ are asymptotic $(r,b(r,k))$-approximate groups 
for all $r \geq 2$.
\et

\begin{proof}
For every positive integer $h$, we have $|h\ast A | = |A | = k$ and,
by Lemma~\ref{approx:lemma:ConvexSumset}, 
\[
hK = h\ast K = \conv(h\ast A).
\]
It follows from Theorem~\ref{AAG:theorem:PolytopeApproxGroup}
that there exists a set $X_{r,h} \subseteq \mcs(h\ast A) \subseteq \mcs(A)$ 
that satisfies~\eqref{AAG:PolytopeAsymptoticApproxGroup-X} 
and~\eqref{AAG:PolytopeAsymptoticApproxGroup},
and so $hK$ is an $(r,b(r,k))$-approximate group.  

Let $K' = (k-1)K$.  Dilating the sets on both sides 
of~\eqref{AAG:PolytopeAsymptoticApproxGroup} by $k-1$, we obtain 
\[
rh K' = rh (k-1)K \subseteq  X_{r,h}+ h (k-1)K = X_{r,h}+ hK'.  
\]
This completes the proof.  
\end{proof}

\bl              \label{AAG:lemma:(h-c)K}
Let  $A = \{a_0, a_1,\ldots, a_{k-1}\}$ be a finite subset 
of a real vector space $V$, with $|A| = k \geq 2$.
Let $\mcs(A)$ be the subsemigroup of $V$ generated by $A$.  
Consider the polytope 
\[
K = \conv(A).
\]
If $c$ is a positive integer and $h \geq ck$,  
then $c\ast A \subseteq \mcs(A)$ and 
\beq                          \label{AAG:(h-c)K}
hK \subseteq c\ast A + (h-c)K.
\eeq
\el

\begin{proof}
For every positive integer $c$, the dilation $c\ast A$ is a subset of 
the semigroup $\mcs(A)$.  Let $h \geq ck$.  
If $\sum_{i=0}^{k-1} \mu_i = h$ with $\mu_i \geq 0$ for all $i$,   
then  $x = \sum_{i=0}^{k-1} \mu_ia_i \in hK$.  We have  
\[
0 < ck\leq h = \sum_{i=0}^{k-1} \mu_i \leq k \max( \{\mu_i:i=0,1,\ldots, k-1\})
\]
and so $\mu_j \geq c$ for some $j \in  \{ 0,1,\ldots, k-1\}$.  
It follows that 
\[
x' =  \sum_{\substack{i=0 \\ i \neq j}}^{k-1} \mu_ia_i + (\mu_j-c)a_j \in (h-c)K
\]
and 
\[
x = ca_j + x' \in c\ast A + (h-c)K.
\]
This completes the proof. 
\end{proof}

For every real number $x$, the \emph{integer part} of $x$ is the unique integer 
$[x]$ such that $[x] \leq x < [x]+1$, and the \emph{fractional part} of $x$  is 
$(x) = x - [x] \in [0,1)$.
The following result is essentially Propositions 1 and 2  
of Khovanskii~\cite{khov92}.    
We follow the proof in~\cite{khov92}.

\bl         \label{approx:lemma:LatticePointRewrite}
Let $A = \{0, a_1,\ldots, a_{k-1}\}$ be a finite subset of a real normed 
vector space $V$ with  $|A| = k \geq 2$ 
such that the subgroup $\mcg(A)$ generated by $A$ 
is a discrete subgroup of $V$.  
Consider the polytope $K = \conv(A)$.  
There exists a positive integer $c = c(A)$ such that, if $h \geq ck$, then   
\beq          \label{approx:LatticePointRewrite}
 \mcg(A) \bigcap \left(  c\sum_{i=1}^{k-1} a_i  +   (h-ck)K \right) 
 \subseteq  hA.
\eeq
\el

\begin{proof}
For all positive integers $c$ and $h$ with $h \geq ck$, we define 
\begin{align*}
K(h,c) & = c\sum_{i=1}^{k-1} a_i  + (h-ck)K \\
& = \left\{  \sum_{i=1}^{k-1} x_ia_i : x_i \geq c \text{ and } 
\sum_{i=1}^{k-1} x_i \leq h-c \right\}.  
\end{align*}
The ``parallelepiped''   
\[
P = \left\{\sum_{i=1}^{k-1} t_i a_i: 0 \leq t_i < 1 \text{ for } i=1,\ldots, k-1 \right\}
\] 
is a bounded subset of $V$,  
and so $P$ contains only finitely many elements of the discrete group $\mcg(A)$.  
Because $A$ generates the  discrete group $\mcg(A)$, 
for each vector $p \in P \cap \mcg(A)$ there exist 
integers $ z_{p,1},\ldots, z_{p,k-1}$ such that $p = \sum_{i=1}^{k-1} z_{p,i} a_i$.
Let 
\[
m = \max\left( \sum_{i=1}^{k-1} |z_{p,i}| : p \in P\cap \mcg(A) \right)
\]
and 
\[
c  = k-1+m.
\]
For all $x_1,\ldots, x_{k-1} \in \R$, we have   
\[
q = \sum_{i=1}^{k-1} x_i a_i 
= \sum_{i=1}^{k-1} [x_i ]a_i + \sum_{i=1}^{k-1} (x_i) a_i 
\]
where 
\[
\sum_{i=1}^{k-1} [x_i ]a_i \in \mcg(A)
\qqand 
\sum_{i=1}^{k-1} (x_i) a_i \in P. 
\]    
If $q \in \mcg(A)$, then 
\[
p = \sum_{i=1}^{k-1} (x_i) a_i = q - \sum_{i=1}^{k-1} [x_i ]a_i 
=  \sum_{i=1}^{k-1} z_{p,i} a_i  \in P \cap \mcg(A)
\]
and so 
\[
q =  \sum_{i=1}^{k-1} [x_i ]a_i +  \sum_{i=1}^{k-1} z_{p,i} a_i 
= \sum_{i=1}^{k-1} \left( [x_i ] + z_{p,i} \right) a_i 
= \sum_{i=1}^{k-1} y_i a_i 
\]
where 
\[
y_i  =  [x_i ] + z_{p,i} \in \Z
\]
for $i= 1,\ldots, k-1$.
Moreover,
\begin{align*}
\sum_{i=1}^{k-1} |x_i - y_i|
& = \sum_{i=1}^{k-1} |x_i - [x_i ] - z_{p,i} | 
 = \sum_{i=1}^{k-1} |  (x_i ) - z_{p,i} | \\
& \leq \sum_{i=1}^{k-1}   (x_i )+ \sum_{i=1}^{k-1}| z_{p,i}| \\
&  < k-1 + m = c.
\end{align*}
If $h \geq ck$ and $q  \in \mcg(A) \cap K(h,c)$, 
then there exist real numbers $x_1,\ldots, x_{k-1}$  such that 
\[
q = \sum_{i=1}^{k-1} x_i a_i 
\]
where $x_i \geq c$ for $i=1,\ldots, k-1$ and $\sum_{i=1}^{k-1} x_i \leq h-c$.  
There also exist integers $y_1,\ldots, y_{k-1}$ such that 
\[
q = \sum_{i=1}^{k-1} y_i a_i
\qqand  
\sum_{i=1}^{k-1} \left| x_i  - y_i \right| < c.
\]
The second inequality implies that  
\[
 x_i  - y_i  \leq \left| x_i  - y_i \right| < c
 \]
and so 
\[
y_i > x_i - c \geq 0.
\]
Thus, the integers $y_1,\ldots, y_{k-1}$  are positive.  
Similarly,  
\[
 \sum_{i=1}^{k-1}  y_i - \sum_{i=1}^{k-1}   x_i 
 = \sum_{i=1}^{k-1}(  y_i - x_i ) \leq  \sum_{i=1}^{k-1} \left| x_i  - y_i \right| < c
\]
and so
\[
h' =  \sum_{i=1}^{k-1}  y_i < \sum_{i=1}^{k-1}   x_i  + c \leq (h-c)+c = h.
\]
Because $0 \in A$ and $h' \leq h$, we have  $h'A \subseteq hA$ 
and 
\[
q =  \sum_{i=1}^{k-1} y_i a_i \in \left( \sum_{i=1}^{k-1}y_i \right)A 
= h'A  \subseteq hA.
\]
Therefore, 
\[
\mcg(A) \cap K(h,c) 
= \mcg(A) \cap \left( c\sum_{i=1}^{k-1} a_i  + (h-ck)K \right)  \subseteq hA.
\]
This completes the proof.  
\end{proof}

\bt             \label{AAG:theorem:Main}
Let $A = \{a_0, a_1,\ldots, a_{k-1}\}$ be a finite subset of a real normed 
vector space $V$ with $|A| = k \geq 2$ 
such that the subgroup $\mcg(A)$ generated by $A$ is a discrete subgroup of $V$. 
Let $A'= (k-1)\ast A$.  
There exists $c = c(A)$ such that, for all $r \geq 2$ and $h \geq ck^2$, 
there is a set $X'_{r,h} \subseteq \mcg(A)$ with $|X'_{r,h}| \leq k b(r,k)$ 
and  
\[
rhA'   \subseteq X'_{r,k}+ hA'.
\]
Thus, the set $(k-1)\ast A$ is an asymptotic $(r,k b(r,k))$-approximate group.  
\et
 
\begin{proof} 
The translated set  $A_0 = A - a_0$
satisfies $|A_0| = k$ and $0 \in A_0$.   
The group $\mcg(A_0)$ is a subgroup of the discrete group $\mcg(A)$, and 
so $\mcg(A_0)$ is also discrete.  
By Lemma~\ref{AAG:lemma:translate}, it suffices to prove that $A_0$ 
is an asymptotic $(r,k b(r,k))$-approximate group.  
Equivalently, we shall assume that $A = \{0,a_1, a_2,\ldots, a_{k-1}\}$, 
and apply Lemma~\ref{approx:lemma:LatticePointRewrite}.  

Let
\[
K = \conv(A)
\]
and let $c = c(A)$ be the integer constructed 
in Lemma~\ref{approx:lemma:LatticePointRewrite}.  
If 
\[
A' = (k-1)\ast A
\]
then $|A'| = |(k-1)\ast A| = k$ and 
\[
K' = \conv(A')  = (k-1)K.
\]
Let  $\mcs(A)$ be the subsemigroup of $V$ generated by $A$.
Applying formula~\eqref{AAG:PolytopeAsymptoticApproxGroup-K'}
in Theorem~\ref{AAG:theorem:PolytopeAsymptoticApproxGroup},  
we obtain, for every integer $h \geq 1$, 
a set $X_{r,h} \subseteq \mcs(A)$ with 
\[
|X_{r,h}| \leq b(r,k)
\]
and
\[
r h K' \subseteq X_{r,h}  + h K'.
\]
Applying formula~\eqref{AAG:(h-c)K} in Lemma~\ref{AAG:lemma:(h-c)K} 
to the polytope $K'$ with $h \geq ck^2$, we obtain 
\[
hK' \subseteq ck \ast  A' + (h-ck)K'
\]
and so
\[
r h K' \subseteq X_{r,h}  + ck\ast  A' + (h-ck)K'.
\]
The finite set 
\[
X'_{r,k} = X_{r,k}  + ck \ast A' - c(k-1) \sum_{i=1}^{k-1} a_i.
\]
is a contained in $\mcg(A')$ and satisfies 
\[
|X'_{r,k}| \leq |A'| |X_{r,k}|  \leq k b(r,k).  
\]
Applying formula~\eqref{approx:LatticePointRewrite} from 
Lemma~\ref{approx:lemma:LatticePointRewrite} 
to the sets $A'$ and $K'$, we obtain  
\begin{align*}
rhA'   \subseteq & \mcg(A') \cap r hK'  \\
 \subseteq &  \mcg(A') \cap \left( X_{r,k}  + ck\ast A' + (h-ck)K' \right) \\
 =  & \left( X_{r,k}  + ck \ast A' - c(k-1) \sum_{i=1}^{k-1} a_i   \right) \\
& +  \mcg(A') \cap \left( c(k-1) \sum_{i=1}^{k-1} a_i   + (h-ck)K' \right) \\
 \subseteq & X'_{r,k}+ hA'
\end{align*}
for all $h \geq ck^2$, and so $A' = (k-1) \ast A$ is an asymptotic $(r, kb(k,r))$-approximate group.   
This completes the proof.  
\end{proof}

\bt             \label{AAG:theorem:Main-Zn}
Let $k \geq 2$ and $r \geq 2$.  
Every set of $k$ lattice points is an asymptotic $(r,k b(k,r))$-approximate group.  
\et

\begin{proof} 
Let $A$ be a subset \Zn\ with $|A| = k$.  
The group \Zn\ is a discrete subgroup of the normed 
vector space $\R^n$, and so the subgroup of \Zn\ generated by $A$
is also discrete.  
Thus, the set $A$ satisfies the conditions 
of Theorem~\ref{AAG:theorem:Main}.     
Let $c = c(A)$ be the positive integer constructed 
in Lemma~\ref{approx:lemma:LatticePointRewrite}, 
and let $A' = (k-1)\ast A$. 
By Theorem~\ref{AAG:theorem:Main}, for $r \geq 2$ and $h \geq ck^2$, 
there is a set $X'_{r,k} \subseteq \Z^n$ 
such that  
\[
 (k-1)\ast  \left( rh A \right) = rh A'   \subseteq X'_{r,k} + hA' = X'_{r,k}+   (k-1)\ast hA.
\]
Lattice points in the sets $ (k-1)\ast  \left( rh A \right) $ and $ (k-1)\ast hA$ 
have all coordinates divisible by $k-1$.   
If $X''_{r,k}$ is the subset of $X'_{r,k}$ consisting only of vectors 
all of whose coordinates are divisible by $k-1$,  then
\[
 (k-1)\ast  \left( rh A \right)    \subseteq X''_{r,k}+   (k-1)\ast hA.
\]
Let $Y_{r,k}$ be the subset of $\Z^n$  
such that $X''_{r,k} = (k-1)\ast Y_{r,k}$.  We have 
\[
|Y_{r,k}| = |X''_{r,k}| \leq |X'_{r,k}| \leq kb(r,k) 
\]
and
 \[
 (k-1)\ast  \left( rh A \right)    \subseteq (k-1)\ast Y_{r,k} +   (k-1)\ast hA.
\]
Dividing by $k-1$, we obtain 
 \[
rh A  \subseteq  Y_{r,k} + hA.
\]
Thus, the set $A$ is an asymptotic $(r,k b(k,r))$-approximate group.  
\end{proof}

\bt
Every nonempty finite subset of an abelian group 
is an asymptotic approximate group.
\et

\begin{proof} 
Consider the group $G_0 \times \Z^n$, 
where $G_0$ is a finite abelian group of order $n_0$, and $n \in \N$.
Let $\pi_1: G_0 \times \Z^n \rightarrow \Z^n$ be the projection homomorphism.  
If $A'$ is a nonempty finite subset of $G_0 \times \Z^n$, then $\pi_1(A')$ 
is a nonempty finite subset of $\Z^n$.  Let $k = |\pi_1(A')|$ and let $r \geq 2$.   
By Theorem~\ref{AAG:theorem:Main-Zn}, the set $\pi_1(A')$ is an 
asymptotic $(r,kb(r,k))$-approximate subgroup of $\Z^n$.  
By Lemma~\ref{AAG:lemma:Main-projection}, the set $A'$ is an 
asymptotic $( r, n_0 k b(r,k))$-approximate subgroup of $G_0 \times \Z^n$.  

Let $A$ be a finite subset of an abelian group, 
and let $G = \mcg(A)$ be the subgroup generated by $A$.
Every finitely generated abelian group is isomorphic to a group 
of the form $G_0 \times \Z^n$, where $G_0$ is a finite group and $n \in \N$.
Let $f:G \rightarrow G_0 \times \Z^n$ be an isomorphism.  
If $A$ is a nonempty finite subset of $G$, then $A' = f(A)$ is a 
nonempty finite subset of $G_0 \times \Z^n$.  
For every integer $r \geq 2$ the set $A'$ is an 
asymptotic $( r, n_0 k b(r,k))$-approximate subgroup of $G_0 \times \Z^n$.
By Lemma~\ref{AAG:lemma:isomorphism}, 
for every integer $r \geq 2$, the set $A$ is also an 
asymptotic $(  r, n_0 k b(r,k))$-approximate subgroup of $G$.   
This completes the proof.  
\end{proof}

\def\cprime{$'$} \def\cprime{$'$} \def\cprime{$'$} \def\cprime{$'$}
\providecommand{\bysame}{\leavevmode\hbox to3em{\hrulefill}\thinspace}
\providecommand{\MR}{\relax\ifhmode\unskip\space\fi MR }
\providecommand{\MRhref}[2]{%
  \href{http://www.ams.org/mathscinet-getitem?mr=#1}{#2}
}
\providecommand{\href}[2]{#2}

\end{document}